\documentclass[a4paper,10pt]{article}
\title{\Large{\textbf{On a new property of primes that leads to a generalization of Cram\'er's conjecture}}}
\author{\large{Nilotpal Kanti Sinha}
\\
\small{Bangalore, India.}}
\begin{document}
\maketitle{\textit{Dedicated to Dr.A.P.J.Abdul Kalam.}}
\section{Introduction}
One of the most famous open problems in the theory of primes is the gap between consecutive primes. Given two consecutive prime numbers $p_{n+1}$ and $p_n$, how large can the difference $ g(n) = p_{n+1} - p_n$ be? The answer to this question has baffled mathematicians for almost a century; however continuous progress has been since 1920 when Harald Cram\'er's proved on the assumption of the Riemann Hypothesis, that 
\begin{equation}
p_{n+1} - p_n = O\Big(p_n^{0.5}\ln p_n\Big)
\end{equation}
(See \cite{hc1}). Currently the best unconditional result is $O\Big(p_n^{0.535}\Big)$ due to R. Baker and G. Harman (See \cite{bh}). 
\\
\\
\textbf{Cram\'er's Conjecture:} Cram\'er applied probabilistic methods to estimate an upper bound on $g(n)$ and in 1936, he suggested that (See \cite{hc2})
\begin{equation}
\limsup_{n \rightarrow \infty}\frac{p_{n+1} - p_n}{\ln ^2 p_n} = 1.
\end{equation}
This statement is known as the Cram\'er's conjecture. Intuitively, this means the gaps between consecutive primes are always small, and it quantifies asymptotically just how small they can be. This conjecture has neither been proved nor disproved. In Cram\'er's model one assumes that the probability that a natural number $x$ is prime is $1/\ln x$. This model was consistent with empirical and  in 1966, Gallagher showed that the Hardy-Littlewood conjecture on $k$-tuples was consistent with Cramer's model (See \cite{hl}). 
\\
\\
\textbf{Cram\'er-Granville's Conjecture:} Despite so much evidence in support of Cram\'er's probabilistic model, in 1985, Helmut Maier proved a result that actually contradicts the predictions of Cram\'er's (See \cite{hm}). The problem with Cram\'er's model is that it fails to take into account divisibility. Thus, for the primes $p > 2$, the probability that $p+1$ is prime is not $1/\ln(p + 1)$as suggested by the Cram´er model but rather 0 since $p+1$ is even. Further $p+2$ is always odd; therefore it is twice as likely to be prime as a random natural number. Thus $n$ and $n+2$ being primes are not independent events.
Based on Maier's result, in 1995, Andrew Granville refined Cram\'er's model and conjectured that
\begin{equation}
\limsup_{n \rightarrow \infty}\frac{p_{n+1} - p_n}{\ln ^2 p_n} \ge 2e^{-\gamma}
\end{equation}
where $\gamma$ is the Euler-Mascheroni constant (See \cite{ag}). The modified statement 
\begin{equation}
p_{n+1} - p_n < M\ln ^2 p_n
\end{equation}
with $M>1$ is known as the Cram\'er-Granville's Conjecture.
\\
\\
\textbf{Firoozbakht's Conjecture:} In 1982, Farideh Firoozbakht made a lesser known but interesting conjecture that 
\begin{equation}
p_n ^{\frac{1}{n}} > p_{n+1} ^{\frac{1}{n+1}}
\end{equation}
for all $n \ge 1$. (See \cite{pr}). If this is indeed true then it can be shown that 
\begin{equation}
p_{n+1} - p_n < \ln ^2 p_n - \ln p_n
\end{equation}
(Lemma 2.3) for all sufficiently large $n$. Firoozbakht claims to have verified this conjecture for all primes up to $10^{12}$. This upper bound is not only stronger that the Cram\'er's conjecture but it also contradicts Granville's limit $2e^{-\gamma}$. Thus either (4) is false or the Firoozbakht's conjecture is false.
\\
\\
\textbf{Scope of this paper:} I present a study of the gap between consecutive primes as a special case of the gap between sequences having a certain property which I call pseudo equidistribution. In section 2, I present two sufficient conditions for Cram\'er's conjecture. In section 3, I introduce the concept of pseudo equidistribution mod 1 and show that prime numbers are a special case of the family of sequence pseudo equidistributed mod . I give two arguments in support of Cram\'er's conjecture and also show that this conjecture can be extended to pseudo equidistributed sequences. The theoretical argument is given in section 4 and the heuristic argument is given in section 5. The heuristic argument in section 5, also supports Firoozbakht's conjecture; and as mentioned above, it also implies that this is an argument against Granville's limit $2e^{-\gamma}$. Finally in section 6, I present a fromal statement of the generalized Cram\'er's conjecture.
\section{Sufficient conditions for Cram\'er's conjecture}
\textbf{Definition 2.1.} \textit{Let $f(n) = p_n ^{\frac{1}{n}}$. If $f(n) > f(n+1)$, we define $a_n$ as
\begin{displaymath}
a_n = -\frac{1}{\ln n}\ln\Bigg(\frac{f(n)}{f(n+1)} - 1\Bigg),
\end{displaymath}
and if $f(n) < f(n+1)$, then $a_n$ is undefined}. 
\\
\\
From the above definition, it follows that if $f(n) > f(n+1)$ then we have
\begin{equation}
p_n ^{\frac{1}{n}} = \Bigg(1 + \frac{1}{n^{a_n}}\Bigg)p_{n+1} ^{\frac{1}{n+1}}.
\end{equation}
\\
We present two equivalent froms of Cram\'er's conjecture depending on weather $f(n) < f(n+1)$ or  $f(n) > f(n+1)$.
\\
\\
\textbf{Lemma 2.2.} \textit {If there exists a positive constant $c_0$ such that
\begin{displaymath}
p_n ^{\frac{1}{n}} > \Bigg(1-\frac{c_0 \ln p_n}{n^2}\Bigg)p_{n+1} ^{\frac{1}{n+1}}
\end{displaymath}
for all sufficiently large $n$ then the Cram\'er'-Granville conjecture is true}.
\\
\\
\textit{Proof}. If the condition mentioned in the lemma is true then  
\begin{displaymath}
p_{n+1} < p_n ^{1+\frac{1}{n}} \Bigg(1-\frac{c_0 \ln p_n}{n^2}\Bigg)^{-(n+1)}
< p_n e^{\frac{\ln n}{n}} \Bigg(1+\frac{c_0\ln p_n}{n} + \frac{c_2\ln^2 p_n}{n^2} \Bigg)
\end{displaymath}
or
\begin{displaymath}
p_{n+1} < p_n \Bigg(1+ \frac{\ln p_n}{n} + \frac{c_1\ln^2 p_n}{n^2}\Bigg) 
\Bigg(1+\frac{c_0\ln p_n}{n} + \frac{c_2\ln^2 p_n}{n^2} \Bigg)
\end{displaymath}
\\
where $c_1$ and $c_2$ are some positive constants. Hence for all sufficiently large $n$,
\begin{equation}
p_{n+1} - p_n < \frac{(1+c_0) p_n\ln p_n}{n} + 1.
\end{equation}
But $\frac{p_n}{n} < (\ln p_n - 1)$ for all sufficiently large $n$. Hence (8) reduces to 
\begin{equation}
p_{n+1} - p_n < (1+c_0)\Big(\ln^2 p_n - \ln p_n\Big) + 1.
\end{equation}
This proves the lemma.
\\
\\
\textbf{Lemma 2.3.} \textit {If for all sufficiently large $n$, $a_n \ge 0$ then Firoozbakht's conjecture $($and hence as a corollary Cram\'er's conjecture$)$ is true.}
\\
\\ 
\textit{Proof}. If $a_n \ge 0$ then
\begin{displaymath}
p_n ^{\frac{1}{n}} = \Bigg(1 + \frac{1}{n^{a_n}}\Bigg)p_{n+1} ^{\frac{1}{n+1}} > p_{n+1} ^{\frac{1}{n+1}}.
\end{displaymath}
Proceeding exactly as in the proof of Lemma 2.2, we obtain
\begin{equation}
p_{n+1} - p_n < \ln^2 p_n - \ln p_n + 1.
\end{equation}
This proves the lemma.
\\
\\
The above give the sufficient (but not necessary) conditions for Cram\'er's conjecture to be true. Also it can be noted from the above lemma that the error term $O(\ln ^2 p_n)$ is optimum and it cannot be lowered. Further since Lemma 2.2 does require the condition $f(n) > f(n+1)$, this lemma may provide us an easier approach to Cram\'er's conjecture as we shall see in section 4.
\section{Primes are pseudo equidistributed}
I had developed the theory given in this section in section 5 of \cite{nw} but for the sake of the completeness of the current paper on its own, I am reproducing the main results in this section.
\\
\\
\textbf{Definition 3.1.}\textit{A sequence of positive real numbers $s_n$ is said to be pseudo equidistributed mod $1$ if $s_n$ is strictly increasing and the sequence of ratios 
\begin{displaymath}
\frac{s_1}{s_n},\frac{s_2}{s_n},...,\frac{s_n}{s_n}
\end{displaymath}
approach unifrom distribution modulo one as} $n \rightarrow \infty$.
\\
\\
\textbf{Theorem 3.2.} \textit {If $s_n$ is a sequence of positive reals such that $\lim_{n \rightarrow \infty}\frac{1}{s_n}=0$ and $\lim_{n \rightarrow \infty}\frac{s_{[nt]}}{s_n} = t$, for every real $t$ in $[0,1]$, then the sequence of ratios 
$\frac{s_r}{s_n}, (r = 1, 2, \ldots, n)$, is pseudo equidistributed mod $1$ and
\\
\begin{equation}
\lim_{n \rightarrow \infty}\frac{1}{n}\sum_{r \le n}f\Big(\frac{s_r}{s_n}\Big) = \int_{0}^{1}f(x)dx.
\end{equation}
where f is any function Reimann integrable in} $[0,1]$.
\\
\\
\textit{Proof}. If the condition mentioned in the statement of the theorem is true then
\\
\begin{displaymath}
\lim_{n \rightarrow \infty}\frac{s_{[bt]} - s_{[at]}}{s_n} = b-a
\end{displaymath}
\\
where $0 \le a < b \le 1$. Therefore as $n \rightarrow \infty$, the probability of finding an integer $r$ such that $a \le \frac{s_r}{s_n} \le b$ approaches  $(b-a)$. Hence the sequence of the ratios $\frac{s_r}{s_n}, r \le n$, is pseudo equidistributed mod 1 and so we can use properties of unifrom distribution mod 1 on the sequence $\frac{s_r}{s_n}$. It is known that if the sequence $b_n$ is unifromly distributed mod 1 and $f$ is Riemann integrable in $[0,1]$ then
\\
\begin{equation}
\lim_{n \rightarrow \infty}\frac{1}{n}\sum_{r \le n}f(b_r) = \int_{0}^{1}f(x)dx
\end{equation}
\\
(See \cite{gr}, Page 3). Replacing $b_r$ by $\frac{s_r}{s_n}$ in (2) we obtain the required result.
\\
\\
\textbf{Corollary 3.3.}\textit{If two sequences are pseudo equidistributed modulo one then their linear combinations are also pseudo equidistributed modulo one}.
\\
\\
\textit{Proof}. Trivial.
\\
\\
Clearly, the sequence of natural numbers is pseudo equidistributed mod 1. Our next theorem shows that primes are also pseudo equidistributed mod 1.
\\
\\
\textbf{Theorem 3.4.}\textit{The sequence of primes is pseudo equidistributed mod $1$}.
\\
\\
\textit{Proof}. Since $p_n \sim n \ln n$ therefore for all $t$, $0<t<1$,
\\
\begin{displaymath}
\lim_{n \rightarrow \infty}\frac{p_{[tn]}}{p_n} = t + \lim_{n \rightarrow \infty}\frac{t \ln t}{\ln n} = t. 
\end{displaymath}
\\
Hence $p_n$ satisfies all the conditions of Theorem 3.2 therefore $p_n$ is is pseudo equidistributed mod 1.
\\
\\
\textbf{Corollary 3.5.} \textit{If $\alpha$ and $\beta$ are constants, not simultaneously zero then},
\\
\begin{equation}
\lim_{n \rightarrow \infty}\frac{1}{n} \sum_{r \le n} f \Big(\frac{\alpha p_r + \beta r}{\alpha p_n + \beta n}\Big)
= \int_{0}^{1}f(x)dx.
\end{equation}
\\
This result gives us a direct relation between the sequence of primes and the sequence of natural numbers and froms the basis of our heuristic argument in support of Cram\'er's conjecture and its generalization.
\\
\\
\textbf{Lemma 3.6.} \textit {If $s_n$ is pseudo equidistributed mod $1$ and $\lim_{n \rightarrow \infty}s_n^{1/n} = 1$ then}
\begin{displaymath}
s_n^{1/n} = \frac{n+1}{n^2}\sum_{r \le n}s_r^{\frac{1}{n}} + o(1).
\end{displaymath}
\\
\textit{Proof}. Taking $f(x)=x^{\frac{1}{n}}$ in Theorem 3.2 and simplifying, we get
\begin{displaymath}
\lim_{n \rightarrow \infty}\frac{1}{n s_n^{1/n}}\sum_{r \le n}s_r^{1/n} = \frac{n}{n+1} = 1.
\end{displaymath}
Hence
\begin{displaymath}
\lim_{n \rightarrow \infty} \Bigg|s_n^{1/n} - \frac{n+1}{n^2}\sum_{r \le n}s_r^{1/n}\Bigg| = 0.
\end{displaymath}
This proves the lemma.
\\
\\
\textbf{Corollary 3.7.}\textit{The following relations hold}.
\begin{equation}
n^{1/n} = \frac{n+1}{n^2}\sum_{r \le n}r^{1/n} + o(1),
\end{equation}
\begin{equation}
p_n^{1/n} = \frac{n+1}{n^2}\sum_{r \le n}p_r^{1/n} + o(1).
\end{equation}
\\
\textit{Proof}. Follows directly from Lemma 3.5 and Lemma 3.6.
\\
\\
The accuracy of the above fromulas can be seen from the following examples. For $n=1048576$, LHS of (14) is 1.00001322082067
where as the summation in the RHS is 1.00001322082781. Similarly for $n=1048576$, LHS of (15) is 1.00001583690296 where as the summation in the RHS is 1.00001576516749. In fact a stronger from of (15) is
\begin{equation}
p_n^{1/n} = \frac{n+1}{n^2}\sum_{r < n}p_r^{1/n} + \frac{1}{n} + \frac{1}{n\ln n} 
+ O\Bigg(\frac{1}{n \ln^2 n}\Bigg).
\end{equation}
However for our subsequent analysis in section 4, (15) is good enough.
\section{A theoretical argument}
Our first argument in support of Cram\'er's conjecture is based on Lemma 2.2. We present a justification to show that the conditions of Lemma 2.2 are likely to be true.
\\
\\
\textbf{Lemma 4.1.}\textit{If $s_n$ is pseudo equidistributed mod $1$ and $\ln s_n = o(n)$ then for all sufficiently large $n$,}
\begin{displaymath}
\frac{h(n)}{h(n+1)} > 1 - \frac{c\ln s_n}{n^2}
\end{displaymath}
where $c>1$ and $h(n) = \frac{n+1}{n^2}\sum_{r \le n}s_r^{\frac{1}{n}}$.
\\
\\
\textit{Proof}. We have
\begin{equation}
\frac{h(n)}{h(n+1)} 
= \frac{(n+1)^3}{n^2 (n+2)}\frac{\sum_{i \le n}s_i^{\frac{1}{n}}}{\sum_{j \le n+1}s_j^{\frac{1}{n+1}}}
= \frac{(n+1)^3}{n^2 (n+2)}\frac{1}{ \frac{\sum_{i \le n}s_i^{\frac{1}{n+1}}}{\sum_{i \le n}s_i^{\frac{1}{n}}}+ \frac{s_{n+1}^{\frac{1}{n+1}}}{\sum_{i \le n}s_i^{\frac{1}{n}}}}.
\end{equation}
\\
Clearly $n < \sum_{r \le n}s_r^{\frac{1}{n+1}} <\sum_{r \le n}s_r^{\frac{1}{n}}$. Since $\ln s_n = o(n)$ therefore there exists a positive constant $c_1$ such that for all sufficiently large $n$,
\\
\begin{displaymath}
s_{n+1}^{\frac{1}{n+1}} < 1 + \frac{\ln s_{n+1}}{n+1} + \frac{c_1 \ln^2 s_{n+1}}{(n+1)^2}
< 1 + \frac{\ln s_{n+1}}{n} + \frac{c_1 \ln^2 s_{n+1}}{n^2}
\end{displaymath}
\\
Hence from (16), for all sufficiently large $n$,
\begin{displaymath}
\frac{h(n)}{h(n+1)} 
> \frac{(n+1)^3}{n^2 (n+2)}\frac{n}{n + s_{n+1}^{\frac{1}{n+1}}}
> \frac{(n+1)^3}{n(n+2)}\frac{1}{n + 1 + \frac{\ln s_{n+1}}{n} + \frac{c_1 \ln^2 s_{n+1}}{n^2}}
\end{displaymath}
\begin{displaymath}
> \frac{(n+1)^2}{n(n+2)}\Bigg(1 + \frac{\ln s_{n+1}}{n^2} + \frac{c_1 \ln^2 s_{n+1}}{n^3}\Bigg)^{-1}
> 1 - \frac{\ln s_{n+1}}{n^2} - \frac{c_2 \ln^2 s_{n+1}}{n^3}.
\end{displaymath}
where $c_2$ is some positive constant. Hence for every constant $c_3 > 1$, there exists a sufficiently large $n$ such that 
\begin{displaymath}
\frac{h(n)}{h(n+1)} > 1 - \frac{c_3\ln s_{n+1}}{n^2}.
\end{displaymath}
Since $s_n$ is pseudo equidistributed mod 1, $\lim_{n \rightarrow \infty} \frac{s_n}{s_{n+1}} = 1$. Also $\ln s_n = o(n)$. Hence we can choose a constant $c_4 > 1$ such that for every sufficiently large $n$, 
\begin{equation}
\frac{h(n)}{h(n+1)} > 1 - \frac{c_4\ln s_n}{n^2}.
\end{equation}
This proves the lemma.
\\
\\
Although $s_n^{1/n} = h(n) + o(1)$, Lemma 4.1 does not directly imply that 
\begin{equation}
\frac{s_n^\frac{1}{n}}{s_{n+1}^\frac{1}{n+1}} > 1 - \frac{c\ln s_n}{n^2}.
\end{equation}
should hold for all sufficiently large $n$. However from Lemma 3.6, we can see that 
\begin{displaymath}
\Big|s_n^\frac{1}{n} - h(n)\Big| > \Big|s_{n+1}^\frac{1}{n+1} - h(n+1)\Big|.
\end{displaymath}
Hence it is highly likely that (19) is also true. If this is indeed the case then proceeding exactly as in the proof of Lemma 2.1, we can show that for all sufficiently large $n$,
\begin{equation}
s_{n+1}-s_n < \frac{(2+\epsilon) s_n \ln s_n}{n} + 1.
\end{equation}
where $0< \epsilon <1$; and since primes are pseudo equidistributed mod 1, Cram\'er's conjecture would follow from Lemma 2.2 and will have
\begin{displaymath}  
p_{n+1} - p_n < (2+\epsilon) (\ln^2 p_n-\ln p_n).
\end{displaymath}
\section{A heuristic argument}
Unlike the previous which was based on theoretically grounds, our second argument in support of Cram\'er's conjecture is based on a heuristic study of the parameter $a_n$ hence this evidence may not be as strong as the first evidence, none the less, this is still an evidence in favor of Firoozbakht's and consequently Cram\'er's conjecture. This it is also evidence against Granville's conjecture. But before we proceed, I would like to mention a word of caution that $heuristics = o(rigor)$ and number theory full of examples where long standing heuristic predictions have been proven false, the most famous of them being Littlewood's proof of the fact that $\pi(x) > Li(x)$ for infinitely many $x$. 
\\
\\
Since both, the sequence of natural numbers and the sequence of primes are pseudo equidistributed mod 1, we expect at least some of their properties to be analogous. This analogy is seen in several examples. For example 
\begin{equation}
\sum_{r \le n}\frac{1}{r} = \ln n + \gamma + O\Bigg(\frac{1}{n}\Bigg)
\end{equation}
where $\gamma$ is the Euler's constant is analogous to
\begin{equation}
\sum_{r \le n}\frac{1}{p_r} = \ln \ln p_n + M + O\Bigg(\frac{1}{\ln p_n}\Bigg)
\end{equation}
where $M$ is the Merten's constant. Similarly from Corollary 3.5, we see that the relation 
\begin{equation}
\sum_{r \le n}{r}^a \sim \frac{n^{a+1}}{a+1}
\end{equation}
where $a \ne -1$, is analogous to
\begin{equation}
\sum_{r \le n}{p_r}^a \sim \frac{n p_n^a}{a+1}.
\end{equation}
It is on this basis of this analogy that we provide our second evidence. Let $g(n) = n^{\frac{1}{n}}$. We define the sequence $b_n$ for all $n \ge 3$ as  
\begin{displaymath}
b_n = -\frac{1}{\ln n}\ln\Bigg(\frac{g(n)}{g(n+1)} - 1\Bigg).
\end{displaymath}
It is known that $g(n)$ is strictly decreasing for all $n \ge 3$ and hence $b_n > 0$ for all $n \ge 3$. Expanding $b_n$ as series, and taking the dominant terms, we obtain
\begin{equation}
b_n \sim 2 - \frac{\ln \ln n}{\ln n}.
\end{equation}
Hence $b_n \rightarrow 2$ as $n\rightarrow \infty$. Calculating $b_i$ for $i \le 1048576$ we observe the following.
\begin{displaymath}
Mean(b_i) \approx 1.80732285747314
\end{displaymath}
\begin{displaymath}
Median(b_i) \approx 1.81025121723487
\end{displaymath}
The above data shows that the mean and the median are roughly equal. I also observed that this statistical regularity is preserved over any sufficiently large intervals. Of course the approximate values of the mean and the median change according to the interval over which the observation is taken. But in each case, the mean and the median are roughly equal.
\\
\\
I perfromed the same experiment on $a_n$ , as defined in Definition 2.1, for the first 1048576 primes and found that the results are analogous.
\begin{displaymath}
Mean(a_i) \approx 1.79186115958409
\end{displaymath}
\begin{displaymath}
Median(a_i) \approx 1.79480436734964
\end{displaymath}
Not only the mean and the median of $a_i$ are roughly equal but they are also close to the mean and the median of $b_i$; i.e. $a_n \approx 2 - \frac{\ln \ln n}{\ln n}$. This heuristic evidence shows that $a_n \rightarrow 2$ as $n \rightarrow \infty$.
This suggests that the Firoozbakht's conjecture should be true while the Cram\'er's-Granville conjecture is false.
\\
\\
Hence if $f(n) = p_n^{1/n}$ then based on our heuristic analysis, we have
\begin{displaymath}
a_n = \frac{-1}{\ln n} \ln \Bigg(\frac{f(n)}{f(n+1)}-1\Bigg) \sim 2 - \frac{\ln\ln n}{\ln n}
\end{displaymath}
or equivalently
\begin{displaymath}
\lim_{n \rightarrow \infty} \frac{n^2}{\ln n} \Bigg(\frac{f(n)}{f(n+1)}-1\Bigg) = 1.
\end{displaymath}
This would imply that for every $\epsilon, 0<\epsilon<1$, there exists a sufficiently large natural number $N_{\epsilon}$, which depends only on $\epsilon$, such that for all $n>N_{\epsilon}$,
\begin{equation}
\Big(1+n^{-2}\Big) p_{n+1} ^{\frac{1}{n+1}} < p_n ^{\frac{1}{n}} < 
\Big(1+n^{-2+\epsilon}\Big)p_{n+1} ^{\frac{1}{n+1}}.
\end{equation}
This is a stronger form of Firoozbakht's conjecture from which we can deduce that   
\begin{equation}
p_{n+1} - p_n < \ln ^2 p_n - 2\ln p_n + 1.
\end{equation}
\section{Pseudo equidistribution conjectures}
On the basis of the arguments presented in favor of Cram\'er's conjecture in section 4 and section 5, we fromulate the following conjectures on pseudo equidistributed sequences. The Weak conjecture is based on the theoretical argument in section 4 while the Strong conjecture is based on the heuristic argument in section 5. Both these conjecture would imply Cram\'er's conjecture but only the Weak conjecture would imply the Cram\'er'-Granville conjecture.
\\
\\
\textbf{Weak Conjecture 6.1.} \textit{If $s_n>n$ is pseudo equidistributed mod $1$ and $\ln s_n = o(n)$ then for all sufficiently large} $n$,
\begin{displaymath}
s_{n+1} - s_n < \frac{(2+\epsilon) s_n \ln s_n}{n}.
\end{displaymath}
\\
\textbf{Strong Conjecture 6.2.} \textit{If $s_n>n$ is pseudo equidistributed mod $1$ and $\ln s_n = o(n)$ and $g(n) = s_n^{1/n}$ then}
\begin{displaymath}
\frac{-1}{\ln n} \ln \Bigg(\frac{g(n)}{g(n+1)}-1\Bigg) \sim 2 - \frac{\ln\ln n}{\ln n}.
\end{displaymath}
From the strong conjecture we can deduce that
\begin{equation}
s_{n+1} - s_n < \frac{s_n}{n}(\ln s_n - 1) + \frac{c s_n \ln^2 s_n}{n^2}.
\end{equation}
where $c$ is a positive constant.
\section{Conclusion}
In this paper, I have given a new perspective of Cram\'er-Granville's conjecture that leads to its generalization. I believe that the easiest way to attack this conjecture is by proving (19). At one point of time, I thought I had a proof of (19) but I discovered a flaw in the reasoning. This forced me to change the subject of the paper from a proof of the generalized Cram\'er's conjecture to the current title. If  this is a minor flaw that can be rectified then we would prove not only the Cram\'er's conjecture but also its generalization, the Pseudo Equidistribution conjecture.
\section{Acknowledgement}
I am grateful to Marek Wolf for proof reading the paper and pointing out some errors.

\small{e-mail: \texttt{nilotpalsinha@gmail.com}}
\end{document}